\documentclass[aop]{imsart}

\RequirePackage{amsthm,amsmath,amsfonts,amssymb}
\RequirePackage[numbers]{natbib}
\RequirePackage[colorlinks,citecolor=blue,urlcolor=blue]{hyperref}
\RequirePackage{graphicx}

\startlocaldefs
\theoremstyle{plain}

\theoremstyle{remark}

 \newcommand{\be}{\begin{equation}}
 \newcommand{\ee}{\end{equation}}

 \newcommand{\E}{\mbox{\rm \hspace*{.2ex}I\hspace{-.5ex}E\hspace*{.2ex}}}

 \newcommand{\R}{\mbox{\rm \hspace*{.2ex}I\hspace{-.5ex}R\hspace*{.2ex}}}

 \newcommand{\Cov}{\mbox{\rm \hspace*{.2ex}Cov\hspace*{.2ex}}}

 \newtheorem{theo}{\small\bf Theorem}[section]
 \newtheorem{lem}{\small\bf Lemma}[section]
 
 \newtheorem{rem}{\small\bf Remark}[section]
 \newenvironment{REM}{\begin{rem} \rm}{\end{rem}}
 \newtheorem{exam}{\small\bf Example}[section]
 
  \newtheorem{defi}{\small\bf Definition}[section]

 \renewcommand{\Pr}{\mbox{\rm  \hspace*{.2ex}I\hspace{-.5ex}P\hspace*{.2ex}}}

 \newenvironment{pr}[1]{{\small\bf {#1}:}}{}

\endlocaldefs

\begin{document}

\begin{frontmatter}
\title{On corrected Poisson approximations for sums
 of independent indicators}
\runtitle{Corrected Poisson approximations}

\begin{aug}
\author[A]{\fnms{Nickos}~\snm{Papadatos}\ead[label=e1]{npapadat@math.uoa.gr}}
\address[A]{National and Kapodistrian
 University of Athens\printead[presep={,\ }]{e1}}

\end{aug}

\begin{abstract}
 Let $S_n=I_1+\cdots+I_n$ be a sum of independent indicators $I_i$, with
 $p_i=\Pr(I_i=1)=1-\Pr(I_i=0)$, $i=1,\ldots,n$. It is well-known
 that the total variation distance between $S_n$ and $Z_\lambda$, where
 $Z_\lambda$ has a Poisson distribution with mean $\lambda=\sum_{i=1}^n p_i$,
 is typically of order $\sum_{i=1}^n p_i^2$.
 In the present work we propose a class of corrected Poisson
 approximations, which enable the second order factorial moment
 distance  (and hence, the total variation distance)
 to be bounded
 above by a constant multiple of $\sum_{i=1}^n p_i^3$ and
 $\sum_{i=1}^n p_i^4$, hence
 improving the order of approximation.
\end{abstract}

\begin{keyword}[class=MSC]
 \kwd[Primary ]{62F15} 
 \kwd{60E05}
\end{keyword}

\begin{keyword}
 \kwd{Corrected Poisson}
 \kwd{Poisson-Binomial Distribution}
 \kwd{Factorial Moment Distance}
 \kwd{Total Variation Distance}
 \kwd{Gini-Kantorovich-Wasserstein distance}
 \kwd{Inequalities for Elementary Symmetric Functions}
 \kwd{Sums of Independent Indicators}
 \kwd{Improved Order of Approximation}
\end{keyword}

\end{frontmatter}

\section{Introduction}
 \label{sec.intro}
 The total variation distance between two real-valued random variables,
 $X_1,X_2$, is defined by
 \[
 d_{{\rm tv}}(X_1,X_2)=\sup_B |\Pr(X_1\in B)-\Pr(X_2\in B)|,
 \]
 where the supremum is taken over all Borel subsets $B$ of
 $(-\infty,\infty)$. In the particular case where $X_1,X_2$ are non-negative
 integer-valued random variables
 with probability mass-functions $f_1(k)$, $f_2(k)$, $k=0,1,\ldots$,
 the formula simplifies to
 \[
 d_{{\rm tv}}(X_1,X_2)=\frac{1}{2}
 \sum_{k=0}^{\infty}
 \Big|f_1(k)-f_2(k)\Big|.
 \]

 Let $\{I_i\}_{i=1}^n$ be a collection of independent $0-1$
 indicators, with $p_i=\Pr(I_i=1)=1-\Pr(I_i=0)$, $i=1,\ldots,n$.
 The distribution of $S_n=\sum_{i=1}^n I_i$ is
 concentrated on $\{0,\ldots,n\}$
 and it is called Poisson-binomial distribution. Clearly, it is a
 generalization of the Binomial distribution, for if the $p_i$'s are
 all equal (say $p_i=p$ for each $i$),
 then $S_n$ follows a ${\rm Bin}(n,p)$ distribution.
 The distribution of $S_n$ is quite complicated in general, and a
 classical result of Poisson roughly states that if $n$ is large and
 the $p_i$'s are small, then $S_n$ is close to $Z_\lambda$, where
 $Z_\lambda$ is a Poisson random variable with mean
 $\lambda=\E S_n=\sum_{i=1}^n
 p_i$. The quality of Poisson approximation is of fundamental interest, and
 has become a classical theme in applied probability. Among the first
 who
 gave explicit bounds was  Khintchin (1933), LeCam (1960) and Deheuvels and Pfeifer (1986),
 who showed that
 \[
  d_{{\rm tv}}(S_n,Z_\lambda)\leq \sum_{i=1}^n p_i^2.
 \]
  Kontoyiannis et al. (2005), using information theoretic arguments,
 showed  that a similar bound for the Hellinger distance
 is valid,
 namely,
 \[
 d_H^2(S_n,Z_{\lambda}):=\frac{1}{2}\sum_{k=0}^\infty
 \left(\sqrt{\Pr(S_n=k)}-\sqrt{\Pr(Z_{\lambda}=k)}\right)^2\leq
 \frac{1}{\lambda}\sum_{i=1}^n \frac{p_i^3}{1-p_i}.
 \]
 Since $d_{\rm tv}\leq d_{H}\sqrt{2-d_H^2}$, see Novak (2019), this inequality
 for equal $p_i=\lambda/n$ yields the bound
 \[
 d_{\rm tv}(S_n,Z_{\lambda})\leq
\frac{\sqrt{\lambda}\sqrt{2n(n-\lambda)-\lambda^3}}{n(n-\lambda)}\sim
\frac{\sqrt{2\lambda}}{n}, \ \ \mbox{as } \ n\to\infty.
 \]

 The celebrated
 Chen-Stein method, introduced by Chen (1975),
 fruitfully applies to this problem, and provides many other
 approximation results for dependent indicators too; see Barbour et al. (1992).
 One of the
 well-known results (in the independent case)
 is the double inequality
 \[
  \frac{\min\{1,1/\lambda\}}{32}\sum_{i=1}^n p_i^2\leq
 d_{{\rm tv}}(S_n,Z_\lambda)\leq \frac{1-e^{-\lambda}}{\lambda}\sum_{i=1}^n p_i^2,
 \]
 due to Barbour and Eagleson (1983) and  Barbour and Hall (1984).
 Hence, the order of Poisson approximation cannot be improved. For
 example, if $p_i=\lambda/n$ for each $i$, so that $S_n$ is
 ${\rm Bin}(n,\lambda/n)$, then  $d_{{\rm tv}}(S_n,Z_\lambda)\sim
 n^{-1}$, in the sense that $n d_{{\rm tv}}(S_n,Z_\lambda)$ is
 bounded away from $0$ and $\infty$ as $n\to\infty$.
 For a comprehensive
 review in various Poisson approximation results the reader is refered to  Novak (2019); see also
 Roos (1999, 2001), Serfling (1975, 1978).

 The purpose of the present work is to introduce a methodology in order to
 improve the rate of
 convergence from, roughly, $n^{-1}$ to $n^{-2}$ (Theorem \ref{theo.main})
 and to $n^{-3}$ (Theorem \ref{theo.main2}).
 Though possible, we were not able to treat
 either higher orders (except for equal $p_i$ --
 Example \ref{exam.binomial.exact}),
 or dependent indicators.
 The improvements are attained by using suitable signed Poisson measures, which
 we term {\it Corrected Poisson Distributions}, and refer
 to a stronger metric, the factorial moment distance;  see Definition \ref{def.d2}. The main results are based on some novel
 accurate inequalities for the factorial moments of $S_n$ (Lemmas \ref{lem.comparison} and \ref{lem.comparison2}); these inequalities,
 in fact, relate
 the well-known {\it elementary symmetric functions} of $(p_1,\ldots,p_n)$
 with their power sums, $\sum_i p_i^k$, when $p_i\in [0,1]$
 for all $i$.

 In section \ref{sec.last} we give a number of remarks and examples;
 in particular, the binomial case is treated in detail in Example
 \ref{exam.binomial.exact}.

 \section{The corrected Poisson distributions}
 \label{sec.2}
 The second order corrected Poisson distribution is defined as
 \[
 g_{\lambda;\gamma}(k)=e^{-\lambda}\frac{\lambda^k}{k!}
 \Big(1-\gamma\big((k-\lambda)^2-k\big)\Big),
 \ \ \ k=0,1,\ldots \ ,
 \]
 where $\gamma\in\R$ and $\lambda>0$ are constants.
 If $\gamma=0$, $g_{\lambda;0}$ reduces to the ordinary Poisson($\lambda$)
 probability mass function, which we simply denote by
 $\phi_1$ or $g_\lambda$. However, for $\gamma>0$ the values $g_{\lambda;\gamma}(k)$
 become negative for large $k$, hence, it is not a proper
 probability mass function in general, although
 $\sum_{k=0}^\infty g_{\lambda;\gamma}(k)=1$ for all $\gamma\in\R$
 and $\lambda>0$.
 Hence, it will be convenient to make use of the class ${\cal F}_2$,
 defined
 below.
 \begin{defi}
 \label{def.F}
 \[
 {\cal F}_2:=\Big\{g:\{0,1,\ldots\}\rightarrow\R: \sum_{k=0}^\infty |g(k)| u^k<\infty
 {\rm \ for \  some \ } u>2 {\rm \ and \ } \sum_{k=0}^\infty g(k)=1\Big\}.
 \]
 \end{defi}
 It is obvious that the class ${\cal F}_2$ includes both $g_{\lambda;\gamma}$ and
 $f_n$, the probability mass function of $S_n$. Moreover, the
 factorial moments of any function $g\in {\cal F}_2$ can be defined
 as
 \[
 \mu_m(g) := \sum_{k=m}^\infty (k)_m g(k),\ \ m=0,1,\ldots,
 \]
 where $(k)_m:=k(k-1)\cdots(k-m+1)=k!/(k-m)!$ is the descending factorial of
 order $m$ of $k$, with the convention $(k)_0=1$; note that
 all these moments are finite, since the radius of convergence is greater than $2$.
 In the usual proper case
 where $g\geq 0$, $\mu_m(g)=\E (X)_m$ where $X$ follows $g$, and thus,
 $\mu_m(g)$ is just the $m$-th factorial moment of $X$. However, we
 avoid to write $\mu_m=\E (X)_m$ when the distribution of $X$ takes
 negative values (i.e., for signed measures).

 The $\nu$-th order corrected Poisson distribution is defined
 by
 \[
 \phi_{\nu}(k):=e^{-\lambda}\frac{\lambda^k}{k!}
 \Big(1-\gamma_2 P_2(k)-\gamma_3 P_3(k)
 -\cdots-\gamma_{2\nu-2}P_{2\nu-2}(k)
 \Big)\in {\cal F}_2,
 \]
 and it obviously provides an extension of
 $g_{\lambda;\gamma_2}=\phi_2$.
 Here $P_m(k)$ are
 the Poisson-Charlier orthogonal polynomials of
 $Z_\lambda\sim g_{\lambda}=\phi_1$, namely,
 \[
 P_m(k)=\sum_{j=0}^m(-1)^{m-j}{m\choose j} \lambda^{m-j} (k)_j,
 \]
 and satisfy (see, e.g., Afendras et al (2011))
 \[
  \E P_m(Z_\lambda) P_\nu(Z_\lambda)=m!\lambda^m \delta_{m\nu} \ \
  \ \mbox{ and } \ \ \
  \E P_m(Z_{\lambda})g(Z_\lambda)=\lambda^m \E \Delta^m g(Z_{\lambda});
 \]
 observe that the second (covariance) identity implies the first
 (orthogonality) relation, and recall that $\Delta g(k)=g(k+1)-g(k)$ denotes
 the forward difference operator and
 $\Delta^m$ its $m$-th iteration.

 Regarding the total variation distance for functions in ${\cal
 F}_2$, the definition is similar and well understood:
 \begin{defi}{\rm
 \label{def.tvd}
 For $g_1,g_2\in{\cal F}_2$,
 \[
 d_{{\rm tv}}(g_1,g_2):=\frac{1}{2}\sum_{k=0}^{\infty}
 \Big|g_1(k)-g_2(k)\Big|.
 \]
 }
 \end{defi}

 However, in our situation, it is quite complicated to deal with the total variation
 distance, while it seems more convenient to work with a factorial moment distance
 of order
 two, defined as follows:
 \begin{defi}
 \label{def.d2}
 {\rm
 For $g_1,g_2\in{\cal F}_2$,
 \[
 d_2(g_1,g_2):=\frac{1}{2}\sum_{m=1}^{\infty}
 \frac{2^{m}}{m!}\Big|\mu_m(g_1)-\mu_m(g_2)\Big|.
 \]
 }
 \end{defi}
 The metric $d_2$ as well as the following theorem have introduced
 for  the proper case ($g\in{\cal F}_2$, $g\geq 0$) by Afendras and Papadatos (2017);
 the result readily extends to ${\cal F}_2$, hence
 the proof  is omitted.
 \begin{theo}
 \label{theo.d2}
 {\rm
 (a) Any $g\in{\cal F}_2$ can be recovered from its factorial moment
 sequence by the inversion formula
 \[
 g(k)=\frac{1}{k!}\sum_{m=k}^\infty
 \frac{(-1)^{m-k}}{(m-k)!}
 \mu_m(g), \ \ \ k=0,1,\ldots \ .
 \]
 (b)
 For $g_1,g_2\in{\cal F}_2$,
 \[
 d_{{\rm tv}}(g_1,g_2)\leq d_2(g_1,g_2).
 \]
 }
 \end{theo}

 \section{Main results}
 \label{sec.3}
 Let $S_n=\sum_{i=1}^n I_i$ be a sum of $0-1$ independent indicators as in the
 introduction, where $I_i$ has success probability $p_i$, $i=1,\ldots,n$.
 We set $\lambda_j=\sum_{i=1}^n p_i^j$, $j=1,2,\ldots$, and in
 particular, $\lambda_1=\lambda=\E S_n=\sum_{i=1}^n p_i$.
 The accuracy of Poisson approximation cannot be smaller (in magnitude)
 than
 $\sum_{i=1}^n p_i^2$ because the Poisson distribution has equal
 mean and variance, in contrast to $S_n$, which always has smaller
 variance than mean. Due to the tuning parameter $\gamma$, the
 second order corrected Poisson can fill this gap.
 A simple calculation shows
 that
 \begin{lem}
 \label{lem.moments.Poisson}
 {\rm
 \[
 \nu_m:=\mu_m(g_{\lambda;\gamma})=\lambda^m-m(m-1)
 \gamma \lambda^{m},
 \
 \ \ m=0,1,\ldots \ .
 \]
 }
 \end{lem}
 It follows that the "variance" of the corrected Poisson equals to
 $\nu_2+\nu_1-\nu_1^2=(\lambda^2-2\gamma\lambda^2)
 +\lambda-\lambda^2=\lambda-2\gamma\lambda^2$,
 while the variance of $S_n$ is $\sum_{i=1}^n
 p_i(1-p_i)=\lambda-\lambda_2$. Equating variances we are led to the choice
 $\gamma=\lambda_2/(2\lambda^2)$, and our first result reads as follows:
 \begin{theo}
 \label{theo.main}
 {\rm
 Under the preceding notations, let $f_n$ be the
 probability mass function of $S_n$, and set
 \[
 \phi_2(k):=g_{\lambda;\lambda_2/(2\lambda^2)}(k)=
 e^{-\lambda}\frac{\lambda^k}{k!}
 \Big(1-\frac{\lambda_2}{2\lambda^2}\big((k-\lambda)^2-k\big)\Big),
 \ \ \ k=0,1,\ldots \ .
 \]
 Then,
 \[
 d_{{\rm tv}}(f_n, \phi_2)
 \leq
 d_2(f_n,\phi_2)\leq
 \left(\frac{4}{3}\lambda_3+\lambda_2^2\right)e^{2\lambda}
 \leq
 \left(\frac{4}{3}+\lambda\right)
 e^{2\lambda}
 \sum_{i=1}^n p_i^3.
 \]
 }
 \end{theo}
 This inequality provides an essential improvement over
 the traditional Poisson approximation rate, $\sum_{i=1}^n p_i^2$; to see
 this it suffices to take equal $p_i=\lambda/n$. In that case,
 $d_{{\rm tv}}({\rm Bin}(n,\lambda/n),\phi_1)\geq A_\lambda n^{-1}$
 for some constant $A_{\lambda}>0$ depending only
 on $\lambda$, while
 Theorem \ref{theo.main} shows that
 $d_{{\rm tv}}({\rm Bin}(n,\lambda/n),\phi_2)\leq B_{\lambda}
 n^{-2}$.

 At this point we note that a second order correction of similar
 nature is provided by Theorem 3 in Barbour and Hall (1984); namely, for
 $A\subseteq\{0,1,\ldots\}$ they defined
 the quantity
 \[
 \Delta(A):=\Pr(S_n\in A)-\Pr(Z_\lambda\in A)+\frac{\lambda_2}{2\lambda^2}
  \E\left\{I_A(Z_{\lambda})\Big(Z_\lambda^2-(2\lambda+1)Z_{\lambda}+
 \lambda^2\Big)\right\}
 \]
 and proved that
 \[
 \sup_A \Big|\Delta(A)\Big| \leq A_{\lambda} \sum_{i=1}^n p_i^3+B_{\lambda} \
 \Big(\sum_{i=1}^n p_i^2\Big)^2\leq (A_\lambda+\lambda B_\lambda)\sum_{i=1}^n p_i^3.
 \]
 In our notation, $\sup_A \big|\Delta(A)\big|
 =d_{\rm tv}(f_n,\phi_2)$.
 Also, the quantity $\sum_{i=1}^n p_i^3$ appears also in Barbour and Hall's correction, and their constant is (much) smaller than the one provided by Theorem \ref{theo.main}; thus, the bound regarding the total
 variation
 distance is not new. However, the result concerning the stronger metric
 $d_2$ is novel, and can be extended to  higher orders,
 see, e.g., Theorem \ref{theo.main2} below.
 To this end, we choose the constants $\gamma_2,\gamma_3,\gamma_4$
 in order to
 fit the moments of $S_n$ up to order three,
 and we obtain the third order corrected Poisson distribution
 \begin{eqnarray*}
 \phi_3(k) & := &
 e^{-\lambda}\frac{\lambda^k}{k!}
 \Big(1
 -\frac{\lambda_2}{2\lambda^2} P_2(k)
 +\frac{\lambda_3}{3\lambda^3} P_3(k)
 +\frac{\lambda_2^2}{8\lambda^4} P_4(k)
 \Big)
 \\
 &=& e^{-\lambda}\frac{\lambda^k}{k!}
 \Big(a_0+a_1 k+a_2 (k)_2+a_3 (k)_3+a_4 (k)_4
 \Big)
 \end{eqnarray*}
 where
 \[
 \mbox{$
 a_0=1-\frac{\lambda_2}{2}-\frac{\lambda_3}{3}+\frac{\lambda_2^2}{8},
  \
 a_1=\frac{2\lambda_2-\lambda_2^2+2\lambda_3}{2\lambda},
  \
 a_2=\frac{3\lambda_2^2-2\lambda_2-4\lambda_3}{4\lambda^2},
  \
 a_3=\frac{2\lambda_3-3\lambda_2^2}{6\lambda^3},
  \
 a_4=\frac{\lambda_2^2}{8\lambda^4}.
 $}
 \]
 It is easy to verify the following
 \begin{lem}
 \label{lem.moments.Poisson2}
 {\rm
 \[
 \mu_m(\phi_3)=
 \lambda^m-\frac{(m)_2}{2}\lambda_2\lambda^{m-2}
 +\frac{(m)_3}{3}\lambda_3\lambda^{m-3}+
 \frac{(m)_4}{8}\lambda_2^2\lambda^{m-4},
 \ \ m=0,1,\ldots \ ,
 \]
 where $(m)_j=m!/(m-j)!$.
 }
 \end{lem}

 Using the above third-order correction we have the following result.
 \begin{theo}
 \label{theo.main2}
 {\rm
 Under the preceding notations and the assumptions
 of Theorem \ref{theo.main},
 \[
 d_{{\rm tv}}(f_n,\phi_3)\leq d_2(f_n,\phi_3)\leq \frac{2}{3}
 (\lambda^2+4\lambda+3) e^{2\lambda}
 \sum_{i=1}^n p_i^4.
 \]
 }
 \end{theo}
 This inequality provides another essential improvement over
 the previous results, as it is easily verified from
 the particular case of equal $p_i=\lambda/n$. In that case,
 Theorem \ref{theo.main2} implies that
 $d_{{\rm tv}}({\rm Bin}(n,\lambda/n),\phi_3)\leq C_{\lambda}
 n^{-3}$.

 \section{Proofs}
 \label{sec.4}

 The proofs of Theorems \ref{theo.main} and
 \ref{theo.main2} use as a main tool the distance $d_2$,
 and require some accurate bounds of the moment sequence
 $\mu_m:=\E (S_n)_m$; these bounds may be of some independent
 interest.
 Firstly we state and prove an auxiliary result.

 \begin{lem}
 \label{lem.positivity}
 {\rm
 For each $j=0,\ldots,m$, $s=1,2,\ldots$, define the quantities
 \[
 \theta_j(m,s):=\sum_{k=j}^m (-1)^{k-j} {m \choose k}
 \lambda_{k+s} \lambda^{m-k}.
 \]
 Then, $\theta_j(m,s)\geq 0$.
 }
 \end{lem}
 \begin{pr}{Proof}
 Consider a random variable $X$ taking values in $\{1,\ldots,n\}$ with
 respective probabilities $\pi_i:=p_i/\lambda$, $i=1,\ldots,n$, and set
 $Y:=h(X)$, where $h(i):=\pi_i$, $i=1,\ldots,n$. Then,
 \[
 \E Y^{k+s-1}=\sum_{i=1}^n \pi_i h(i)^{k+s-1}
 =\sum_{i=1}^n \pi_i^{k+s}
 =\frac{\lambda_{k+s}}{\lambda^{k+s}}.
 \]
 We thus have
 \[
 \frac{ \theta_j(m,s)}{\lambda^{m+s}}=
 \E\left\{Y^{s-1}
 \sum_{k=j}^m (-1)^{k-j} {m \choose k}
 Y^k\right\}.
 \]
 This is clearly nonnegative when $j=0$, because it represents the expectation
 of $Y^{s-1}(1-Y)^m$ and $0\leq Y\leq 1$ by definition.
 For $j\in\{1,\ldots,m\}$ we shall make use of the identity
 \[
 w(x):=\sum_ {k=j}^m (-1)^{k-j} {m \choose k}
 x^k=j{m\choose j} \int_0^x t^{j-1}(1-x+t)^{m-j} dt.
 \]
 This identity can be proved as follows: Let $U_{1:m}<\ldots<U_{m:m}$ be the order statistics
 from the uniform distribution over the interval $(0,1)$.
 For every $a\in(0,1)$,
 \[
 \Pr(U_{j:m}\leq a)=\Pr( {\rm at \ least\  }  j \  {\rm \ out \ of\ } m \
  {\rm  are\ less\ than\ } \ a)
 =\sum_{k=j}^m {m \choose k } a^k (1-a)^{m-k}.
 \]
 Since  it is well-known that
 $U_{j:m}$ follows a Beta$(j,m+1-j)$ density, we also have
 \[
  \Pr(U_{j:m}\leq a)=j{m\choose j}\int_0^a y^{j-1}(1-y)^{m-j} dy,
 \]
 and equating the above expressions we get
 \[
 \sum_{k=j}^m {m \choose k } a^k (1-a)^{m-k}
 =j{m\choose j}\int_0^a y^{j-1}(1-y)^{m-j} dy.
 \]
 Because both sides represent entire functions of $a$ in the complex plan
 (polynomials), the identity holds
 for all $a$, and we are allowed to set $a=-x/(1-x)$, $x\neq 1$,
 obtaining
  \[
 \sum_{k=j}^m (-1)^k {m \choose k }  x^k
 =j{m\choose j}(1-x)^m \int_0^{-x/(1-x)} y^{j-1}(1-y)^{m-j} dy.
 \]
 The substitution $y=-t/(1-x)$ in the last integral yields
 \[
 \sum_{k=j}^m (-1)^k {m \choose k }  x^k
 =(-1)^j j{m\choose j}\int_0^{x} t^{j-1}(1-x+t)^{m-j} dt,
 \]
 which is equivalent to the desired identity.

 Now, the integral expansion shows that $w(x)\geq 0$ for $0\leq x\leq 1$, and
 we arrived at the representation
 \[
 \frac{ \theta_j(m,s)}{\lambda^{m+s}}=
 \E\left\{Y^{s-1} w(Y) \right\}\geq 0,
 \]
 completing the proof.
 \vspace{1em}
 \end{pr}

 We are now in a position to compare the factorial moments of $S_n$ with those of the
 corrected Poisson  $\phi_2$. Surprisingsly enough, it turns out
 that the
 sequence
 $\mu_m=\E (S_n)_m$
 dominates $\nu_m:=\mu_m(\phi_2)$ for all $m$.
 \begin{lem}
 \label{lem.comparison}
 {\rm
 If $\mu_m$ is the $m$-th factorial moment of $S_n$ then
 \[
 \mu_m=m! S_{n,m},
 \]
 where $S_{n,m}=\sum p_{i_1} \cdots p_{i_m}$ (known as elementary symmetric function of order $m$ in the variables $p_1,\ldots,p_n$) -- the sum runs over all
 ${n \choose m}$ combinations $\{i_1,\ldots,i_m\}$ of $\{1,\ldots,n\}$ (with the convention
 $S_{n,m}=0$ for $m>n$).
 Moreover, the following inequality holds true for all $m\geq 1$:
 \[
  \lambda^m-\frac{(m)_2}{2} \lambda_2 \lambda^{m-2}\leq
 \mu_m\leq
 \lambda^m-\frac{(m)_2}{2} \lambda_2 \lambda^{m-2} +
 \frac{(m)_3}{3}
 \lambda_3 \lambda^{m-3}
 +\frac{(m)_4}{8}
 \lambda_2^2 \lambda^{m-4}=\mu_m(\phi_3),
 \]
 where $(m)_2=m(m-1)$, $(m)_3=m(m-1)(m-2)$,  $(m)_4=m(m-1)(m-2)(m-3)$.
 }
 \end{lem}
 \begin{pr}{Proof}
 The expression for $\mu_m$ is well-known, see e.g.\
 Galambos (1987), so it remains to show the inequalities.
 The proof will be done by induction on $m$. Regarding the lower bound, this is obviously
 true for $m=1,2$. Observe that
 \[
 \mu_{m+1}=(m+1)! S_{n,m+1}=\sum_{i=1}^n p_i \Big\{m! S_{n-1,m}(i)\Big\},
 \]
 where $S_{n-1,m}(i)=\sum p_{i_1} \cdots p_{i_m}$, and where now the sum runs over all
 ${n-1 \choose m}$ combinations $\{i_1,\ldots,i_m\}$ of $\{1,\ldots,n\}\setminus\{i\}$.
  Hence, assuming that the lower bound holds for some $m$, we obtain
 \[
 \mu_{m+1}
 \geq
 \sum_{i=1}^n p_i \left\{(\lambda-p_i)^m-\frac{(m)_2}{2}
 (\lambda_2-p_i^2) (\lambda-p_i)^{m-2}\right\}=S, \ {\rm say},
 \]
 and it suffices to show that $S$ is bounded below by
 $\lambda^{m+1}-\frac{(m+1)_2}{2} \lambda_2 \lambda^{m-1}$.
 Expanding the binomial terms and interchanging the order of summation
 we find
 \begin{eqnarray*}
  S &= &\sum_{k=0}^m (-1)^k {m \choose k} \lambda^{m-k}
 \sum_{i=1}^n p_i^{k+1}
 -\frac{(m)_2}{2} \lambda_2 \sum_{k=0}^{m-2} (-1)^k {m-2 \choose k} \lambda^{m-2-k}
 \sum_{i=1}^n p_i^{k+1} \\
 & &
  +\frac{(m)_2}{2}  \sum_{k=0}^{m-2} (-1)^k {m-2 \choose k} \lambda^{m-2-k}
 \sum_{i=1}^n p_i^{k+3}\\
 &=&
 \sum_{k=0}^m (-1)^k {m \choose k} \lambda^{m-k}
 \lambda_{k+1}
 -\frac{(m)_2}{2} \lambda_2 \sum_{k=0}^{m-2} (-1)^k {m-2 \choose k} \lambda^{m-2-k}
 \lambda_{k+1} \\
 & &
  +\frac{(m)_2}{2}  \sum_{k=0}^{m-2} (-1)^k {m-2 \choose k} \lambda^{m-2-k}
 \lambda_{k+3}=S_1+S_2+S_3, \ \ {\rm say}.
 \end{eqnarray*}
 According to Lemma \ref{lem.positivity},
 $S_3=\frac{(m)_2}{2} \theta_0(m-2,3)\geq 0$. Also,
 \[
 S_1= \lambda^{m+1}-m\lambda_2\lambda^{m-1} +\theta_2(m,2),
 \ \ \
 S_2=-\frac{(m)_2}{2} \lambda_2  \left\{\lambda^{m-1}-
  \theta_1(m-2,1)\right\}.
 \]
 Combining the above and using the fact that the $\theta$'s are nonnegative,
 we obtain
 \[
 S\geq  \lambda^{m+1}-m\lambda_2\lambda^{m-1} -\frac{(m)_2}{2}
 \lambda_2  \lambda^{m-1}
 = \lambda^{m+1}-\frac{(m+1)_2}{2}
 \lambda_2  \lambda^{m-1},
 \]
 and the lower bound is proved.
 Regarding the upper bound, the inequality
 $\mu_m\leq\mu_m(\phi_3)$
 is trivially
 true for $m=1,2,3$ (as equality), and for $m=4$
 yields $\mu_4(\phi_3)-\mu_4=6\lambda_4>0$.
 Assuming that the bound holds for some $m\geq 4$ and proceeding as before, we obtain
 \begin{eqnarray*}
 \mu_{m+1}
 & \leq &
 \sum_{i=1}^n p_i \left\{(\lambda-p_i)^m-\frac{(m)_2}{2}
 (\lambda_2-p_i^2) (\lambda-p_i)^{m-2}
 +\frac{(m)_3}{3}
 (\lambda_3-p_i^3) (\lambda-p_i)^{m-3}
 \right.
 \\
 &&
 \hspace*{30ex}
 \left.
 +\frac{(m)_4}{8}
 (\lambda_2-p_i^2)^2 (\lambda-p_i)^{m-4}
 \right\}=S, \ {\rm say},
 \end{eqnarray*}
 and we shall now show that $S$ is bounded above by
 \[
 \lambda^{m+1}-\frac{(m+1)_2}{2} \lambda_2 \lambda^{m-1}
 +\frac{(m+1)_3}{3} \lambda_3 \lambda^{m-2}
 +
 \frac{(m+1)_4}{8} \lambda_2^2 \lambda^{m-3}=\mu_{m+1}(\phi_3).
 \]
 On expanding the binomial terms and by interchanging the order of summation
 we find $S=\sum_{i=1}^6 S_i$
 where
 \begin{eqnarray*}
  S_1 & = &
  \sum_{k=0}^m (-1)^k {m \choose k} \lambda^{m-k}\lambda_{k+1},
  \\
  S_2  & = &
   -\frac{(m)_2}{2} \lambda_2
  \sum_{k=0}^{m-2} (-1)^k {m-2 \choose k} \lambda^{m-2-k} \lambda_{k+1},
  \\
  S_3  & = &
  \frac{(m)_2}{2}
  \sum_{k=0}^{m-2} (-1)^k {m-2 \choose k} \lambda^{m-2-k} \lambda_{k+3},
  \\
  S_4  & = &
  \frac{(m)_3}{3} \lambda_3
  \sum_{k=0}^{m-3} (-1)^k {m-3 \choose k} \lambda^{m-3-k} \lambda_{k+1},
  \\
  S_5  & = &
  -\frac{(m)_3}{3}
  \sum_{k=0}^{m-3} (-1)^k {m-3 \choose k} \lambda^{m-3-k} \lambda_{k+4},
  \\
  S_6  & = &
  \frac{(m)_4}{8}
  \sum_{k=0}^{m-4} (-1)^k {m-4 \choose k} \lambda^{m-4-k} \Big(\lambda_2^2\lambda_{k+1}-2\lambda_2\lambda_{k+3}+\lambda_{k+5}\Big)
  =\frac{(m)_4}{8} T, \ \ \mbox{say}.
  \end{eqnarray*}
  We shall prove that $T\leq \lambda_2^2 \lambda^{m-3}$.
  Indeed, by considering
  the random variable $Y\in[0,1]$ as in the proof of Lemma \ref{lem.positivity}, we have $\E Y^{k}=\lambda_{k+1}/\lambda^{k+1}$, and hence, after some algebra, $T=\lambda^{m+1}
  \E\left\{ (1-Y)^{m-4}(Y^2-\mu)^2\right\}=\lambda^{m+1}a_m$, say,
  where $\mu=\E Y$. On the other hand,
  $\lambda_2^2 \lambda^{m-3}=\lambda^{m+1}\mu^2$, and the desired
  inequality reduces to $a_m\leq \mu^2$;
  however, the sequence $a_m$ is positive decreasing, and it suffices to
  show that $a_4\leq \mu^2$, i.e., $\E Y^4\leq 2\mu\E Y^2$. Since,
  obviously, $\lambda_5\leq 2\lambda_2\lambda_3$ and
  $\E Y^4=\lambda_5/\lambda^5$, $\mu=\lambda_2/\lambda^2$,
  $\E Y^2=\lambda_3/\lambda^3$, the desired inequality is proved.
  Finally, keeping only the first few terms from the $S_i$'s,
  we conclude, in view of Lemma \ref{lem.positivity},
  the inequalities
  \begin{eqnarray*}
  S_1  & = &
  \mbox{$
  \lambda^{m+1}-m\lambda_2\lambda^{m-1}+\frac{(m)_2}{2}\lambda_3\lambda^{m-2}
  -\theta_3(m,1)
  \leq\lambda^{m+1}
  -m\lambda_2\lambda^{m-1}+\frac{(m)_2}{2}\lambda_3\lambda^{m-2},
  $}
  \\
  S_2  & = &
  \mbox{$
   -\frac{(m)_2}{2} \lambda_2
  \left(
  \lambda^{m-1}-(m-2)\lambda_2\lambda^{m-3}+\theta_2(m-2,1)\right)
  \leq
  -\frac{(m)_2}{2} \lambda_2
  \lambda^{m-1}+\frac{(m)_3}{2}\lambda_2^2\lambda^{m-3},
  $}
  \\
  S_3  & = &
  \mbox{$
  \frac{(m)_2}{2}
  \left(\lambda_3\lambda^{m-2}-\theta_1(m-2,3)\right)\leq
  \frac{(m)_2}{2}
  \lambda_3\lambda^{m-2},
  $}
  \\
  S_4  & = &
  \mbox{$
  \frac{(m)_3}{3} \lambda_3
  \left(\lambda^{m-2}-\theta_1(m-3,1)\right)
  \leq
  \frac{(m)_3}{3} \lambda_3
  \lambda^{m-2},
  $}
  \\
  S_5  & = &
  \mbox{$
  -\frac{(m)_3}{3}
  \theta_0(m-3,4)\leq 0,
  $}
  \\
  S_6  & = &
  \mbox{$
  \frac{(m)_4}{8} T\leq
  \frac{(m)_4}{8}
  \lambda_2^2 \lambda^{m-3}.
  $}
  \end{eqnarray*}
  Hence, in view of the relations
  \[
  \mbox{$
  m+\frac{(m)_2}{2}=\frac{(m+1)_2}{2}, \ \ \
   \frac{(m)_2}{2}+\frac{(m)_2}{2}+\frac{(m)_3}{3}=\frac{(m+1)_3}{2},
   \ \ \
   \frac{(m)_3}{2}+\frac{(m)_4}{8}=\frac{(m+1)_4}{8},
   $}
  \]
  the inductional step is completed and the lemma is proved.
  \vspace{1em}
  \end{pr}

 \noindent
 \begin{pr}{Proof of Theorem \ref{theo.main}} From Lemma \ref{lem.moments.Poisson}
 with $\gamma=\lambda_2/(2\lambda^2)$ we get
 \[
 \mu_m(\phi_2)=\lambda^m-\frac{(m)_2}{2}
 \lambda_2\lambda^{m-2},
 \]
 and Lemma \ref{lem.comparison} shows that
 \[
 0\leq \E(S_n)_m-\mu_m(\phi_2)\leq
 \frac{(m)_3}{3} \lambda_3\lambda^{m-3}
 +\frac{(m)_4}{8} \lambda_2^2\lambda^{m-4}, \  \
 m=1,2,\ldots \ .
 \]
 Thus,
 \begin{eqnarray*}
 d_2(f_n,\phi_2)
  &=&
 \frac{1}{2}\sum_{m=1}^\infty
 \frac{2^m}{m!}\left\{
 \E(S_n)_m-\mu_m(\phi_2)\right\}
 \\
 &&
 \leq \frac{\lambda_3}{6}
 \sum_{m=1}^{\infty}\frac{ 2^m (m)_3\lambda^{m-3}}{m!}
 +\frac{\lambda_2^2}{16}
 \sum_{m=1}^{\infty}\frac{ 2^m (m)_4\lambda^{m-4}}{m!}
 =\Big(\frac{4}{3}\lambda_3+\lambda_2^2\Big)e^{2\lambda}.
 \end{eqnarray*}
 In view of Theorem \ref{theo.d2}(b) and the fact that
 $\lambda_2^2\leq \lambda\lambda_3$, the proof is complete.
 \vspace{1em}
 \end{pr}

 One more accurate bound is needed for the proof of
 Theorem \ref{theo.main2}.

 \begin{lem}
 \label{lem.comparison2}
 {\rm
 If $\mu_m$ is the $m$-th factorial moment of $S_n$ then
 the following lower bound holds true for all $m\geq 1$:
 \begin{eqnarray*}
 \mu_m
  &\geq&
  \lambda^m-\frac{(m)_2}{2} \lambda_2 \lambda^{m-2}
 +\frac{(m)_3}{3} \lambda_3 \lambda^{m-3}
 +
 \frac{(m)_4}{8} \lambda_2^2 \lambda^{m-4}
 -\frac{(m)_4(m)_2}{48} \lambda_4 \lambda^{m-4}.
 \end{eqnarray*}
 }
 \end{lem}
 \begin{pr}{Proof}
 The proof will be done by induction on $m$.
 Denoting  by $L_m$ the lower bound, it is easily checked that
 the inequality $\mu_m-L_m\geq 0$ is satisfied for $m=1,2,3,4$ as an equality. For $m=5$,
 $\mu_5-L_5=4(6\lambda_5+5\lambda\lambda_4-5\lambda_2\lambda_3)$
 and this quantity is positive
 because $\lambda\lambda_4-\lambda_2\lambda_3=\lambda^5 \Cov(Y,Y^2)\geq 0$,
 since $Y\geq 0$ and, thus, both $Y$ and $Y^2$ are increasing functions
 of $Y$ (recall that the random variable $Y\in[0,1]$ is defined in the
 proof of Lemma \ref{lem.positivity}).
 Assuming that the bound is true for some $m\geq 5$ we get
 \begin{eqnarray*}
 \mu_{m+1}& \geq&
 \sum_{i=1}^n p_i \left\{(\lambda-p_i)^m-\frac{(m)_2}{2}
 (\lambda_2-p_i^2) (\lambda-p_i)^{m-2}
 +\frac{(m)_3}{3}
 (\lambda_3-p_i^3) (\lambda-p_i)^{m-3}
 \right.
 \\
 &&\hspace*{1ex}
 \left.
 +\frac{(m)_4}{8}
 (\lambda_2-p_i^2)^2 (\lambda-p_i)^{m-4}
 -\frac{(m)_2 (m)_4}{48}
 (\lambda_4-p_i^4) (\lambda-p_i)^{m-4}
 \right\}
  =S, \ {\rm say}.
 \end{eqnarray*}
 On expanding the binomial terms and interchanging the order of summation
 we can express $S$ as $\sum_{j=0}^{9} S_j$ where
 \begin{eqnarray*}
 S_1& =& \mbox{$
 \lambda^{m+1}-m\lambda_2\lambda^{m-1}
 +\frac{(m)_2}{2}\lambda_3\lambda^{m-2}
 -\frac{(m)_3}{6}\lambda_4\lambda^{m-3}+\theta_4(m,1)$},
 \\
 S_2& =&\mbox{$
  -\frac{(m)_2}{2}\lambda_2\left(\lambda^{m-1}-
 (m-2)
 \lambda_2 \lambda^{m-3}+\frac{(m-2)_2}{2}
 \lambda_3 \lambda^{m-4}-\theta_3(m-2,1) \right)$},
 \\
 S_3& =&\mbox{$
 \frac{(m)_2}{2}
 \left(\lambda_3 \lambda^{m-2}-(m-2)\lambda_4
 \lambda^{m-3}+\theta_2(m-2,3)\right)$},
 \\
 S_4& =&\mbox{$
 \frac{(m)_3}{3}\lambda_3
 \left( \lambda^{m-2}-(m-3)\lambda_2\lambda^{m-4}
 +\theta_2(m-3,1)\right)$},
 \\
 S_5& =&\mbox{$
 -\frac{(m)_3}{3}\left(\lambda_4\lambda^{m-3}-\theta_1(m-3,4)\right),
 $}
  \\
 S_6& =&\mbox{$
 \frac{(m)_4}{8}\lambda_2^2
 \left(\lambda^{m-3}- (m-4)\lambda_2\lambda^{m-5}
 +\theta_2(m-4,1)\right),
 $}
  \\
 S_7& =&\mbox{$
 -\frac{(m)_4}{8}2\lambda_2
 \left(\lambda_3\lambda^{m-4}-\theta_1(m-4,3)\right),
 $}
 \\
 S_8& =&\mbox{$
 \frac{(m)_4}{8}
 \theta_0(m-4,5),
 $}
 \\
 S_9& =&\mbox{$
 -\frac{(m)_4(m)_2}{48}\lambda_4
 \left(\lambda^{m-3}-\theta_1(m-4,1)\right),
  $}
  \\
 S_{0}& =&\mbox{$
 \frac{(m)_4(m)_2}{48}
 \theta_0(m-4,5).
  $}
 \end{eqnarray*}
 Hence, due to the nonnegativity of $\theta$'s, we conclude that
 \begin{eqnarray*}
 S&\geq& \mbox{$
 \lambda^{m+1}-\frac{(m+1)_2}{2}\lambda_2\lambda^{m-1}
 +\frac{(m+1)_3}{3}\lambda_3\lambda^{m-2}
 +\frac{(m+1)_4}{8}\lambda_2^2\lambda^{m-3}$}
\\&&
 \mbox{$
 -\left((m)_3+\frac{(m)_4(m)_2}{48}\right) \lambda_4\lambda^{m-3}
 -\frac{5(m)_4}{6}\lambda_2\lambda_3\lambda^{m-4}
 -\frac{(m)_5}{8} \lambda_2^3\lambda^{m-5}.$}
 \end{eqnarray*}
 Therefore, the inductional step will be proved
 if for all $m\geq 5$,
 \begin{eqnarray*}
 \mbox{$
 \small
 -\left((m)_3+\frac{(m)_4(m)_2}{48}\right) \lambda_4\lambda^{m-3}
 -\frac{5(m)_4}{6}\lambda_2\lambda_3\lambda^{m-4}
 -\frac{(m)_5}{8} \lambda_2^3\lambda^{m-5}
 $}
 &&
 \\
 \mbox{$
  \geq
 -\frac{(m+1)_4(m+1)_2}{48} \lambda_4\lambda^{m-3}.
 $}
 &&
 \end{eqnarray*}
 Multiplying both sides by $\frac{-48}{(m)_4\lambda^{m-5}}$ and collecting terms, we
 arrive at the equivalent inequality ($m=4,5,\ldots$)
 \[
 20\lambda_2\lambda_3\lambda+3(m-4)\lambda_2^3\leq
 (3m+8)\lambda_4\lambda^2.
 \]
 From Cauchy's inequality, $\lambda_2^2\leq
 \lambda_3\lambda$, and hence,
 $\lambda_2^3\leq \lambda_2\lambda_3\lambda$. Therefore,
 \[
 20\lambda_2\lambda_3\lambda+3(m-4)\lambda_2^3\leq
 (3m+8)\lambda_2\lambda_3\lambda,
 \]
 and it suffices to show that $\lambda_2\lambda_3\lambda\leq \lambda_4\lambda^2$, i.e.,
 $\lambda_2\lambda_3\leq \lambda\lambda_4$.
 But this reduces to $\Cov(Y,Y^2)\geq 0$, which is certainly true, and
 the lemma is proved.
 \vspace{1em}
 \end{pr}

 \noindent
 \begin{pr}{Proof of Theorem \ref{theo.main2}} From Lemma \ref{lem.moments.Poisson2}
 we have
 \[
 \mu_m(\phi_{3})=\lambda^m-\frac{(m)_2}{2}
 \lambda_2\lambda^{m-2}+\frac{(m)_3}{3}
 \lambda_3\lambda^{m-3}+\frac{(m)_4}{8}
 \lambda_2^2\lambda^{m-4},
 \]
 and Lemmas \ref{lem.comparison}, \ref{lem.comparison2}, show that
 \[
 0\leq \mu_m(\phi_{3})-\E(S_n)_m\leq \frac{(m)_4(m)_2}{48}
 \lambda_4\lambda^{m-4}, \  \
 m=1,2,\ldots \ .
 \]
 Thus,
 \[
 d_2(f_n,\phi_{3})=\frac{1}{2}\sum_{m=1}^\infty
 \frac{2^m}{m!}\left\{\mu_m(\phi_3)-\E(S_n)_m\right\}
 \leq \frac{1}{2}\lambda_4\sum_{m=0}^{\infty}\frac{2^{m+4}
 \lambda^{m}}{m!}\frac{(m+4)_2}{48},
 \]
 and in view of Theorem \ref{theo.d2}(b), the proof is complete.
 \vspace{1em}
 \end{pr}

 \section{Concluding remarks and examples}
 \label{sec.last}

 \begin{exam}
 \label{exam.binomial}
 {\rm
 Assume that $n>\lambda$.
 The original upper bound for the total variation distance between the Poisson$(\lambda)$
 and Bin$(n,\lambda/n)$  reads as
 \[
 \mbox{$
  \frac{1}{2}\sum_{k=0}^\infty\left|{n\choose k}
 \Big(\frac{\lambda}{n}\Big)^k\Big(1-\frac{\lambda}{n}\Big)^{n-k} -e^{-\lambda}
 \frac{\lambda^k}{k!}
 \right|\leq \lambda(1-e^{-\lambda})n^{-1},
 $}
 \]
 the second order approximation of Theorem \ref{theo.main} implies the bound
 \[
 \mbox{$
  \frac{1}{2}\sum_{k=0}^\infty\left|{n\choose k}
 \Big(\frac{\lambda}{n}\Big)^k\Big(1-\frac{\lambda}{n}\Big)^{n-k} -e^{-\lambda}
 \frac{\lambda^k}{k!}
 \Big( 1-\frac{(k-\lambda)^2-k}{2n}\Big)
 \right|\leq \lambda^3\Big(\frac{4}{3}+\lambda\Big)e^{2\lambda}n^{-2},
 $}
 \]
 while the third order corrected Poisson approximation of
 Theorem \ref{theo.main2} yields
 the bound (here $\lambda_2=\lambda^2/n$,
 $\lambda_3=\lambda^3/n^2$, $\lambda_4=\lambda^4/n^3$)
 \begin{eqnarray*}
 \mbox{$
  \frac{1}{2}\sum_{k=0}^\infty\left|{n\choose k}
 \Big(\frac{\lambda}{n}\Big)^k\Big(1-\frac{\lambda}{n}\Big)^{n-k} -e^{-\lambda}
 \frac{\lambda^k}{k!}\Big( 1-\frac{(k-\lambda)^2-k}{2n}
 +\frac{a_0+a_1 k+a_2 (k)_2+a_3 (k)_3+a_4 (k)_4}{24 n^2}\Big)
 \right|
 $}
 \\
  \mbox{$
 \leq \frac{2}{3}\lambda^4(\lambda^2+4\lambda+3)e^{2\lambda} n^{-3},
 $}
 \end{eqnarray*}
  where  $a_0=\lambda^3(3\lambda-8)$, $a_1=12\lambda^2(\lambda-2)$,
  $a_2=6\lambda(3\lambda-4)$, $a_3=4(2-3\lambda)$, $a_4=3$.
 }
 \end{exam}

 \begin{REM}
 \label{rem.g3tilde}
 {\rm
 Instead of $\phi_3$ of Theorem \ref{theo.main2}, it may appear
 more natural to consider a simpler
 version of the corrected Poisson distribution of order three, namely,
 \[
  \widetilde{\phi}_3(k)=e^{-\lambda}\frac{\lambda^k}{k!}
  \left(1-\frac{\lambda_2}{2\lambda^2}P_2(k)
  +\frac{\lambda_3}{3\lambda^3}P_3(k)\right).
 \]
 Its moments are given
 by $\mu_m(\widetilde{\phi}_3)=\lambda^m-\frac{(m)_2}{2}\lambda_2
 \lambda^{m-2}
 +\frac{(m)_3}{3}\lambda_3\lambda^{m-3}$, and thus,
 $\mu_m=\mu_m(\widetilde{\phi}_3)$, for $m=0,1,2,3$.
 Moreover, the function $\widetilde{\phi}_3$
 is the unique member of the parametric class
 \[
 e^{-\lambda}\frac{\lambda^k}{k!}
  \Big(1-\gamma_1 P_1(k)-\gamma_2 P_2(k)-\gamma_3 P_3(k)\Big)
 \]
 with moments equal to those of $S_n$ up to order three, and one
 might
 expect a third order approximation, comparable to
 that of Theorem \ref{theo.main2} for $\phi_3$. However,
 surprisingly enough, the approximation that $\widetilde{\phi}_3$
 attains is only of order $n^{-2}$, and the adjunction of the polynomial
 $P_3$ does not seem helpful. To see this, it suffices to calculate
 (in the case of equal $p_i=\lambda/n$) the difference
 \[
 f_n(0)-\widetilde{\phi}_3(0)
 =\Big(1-\frac{\lambda}{n}\Big)^n-e^{-\lambda}
 \Big(
 1-\frac{\lambda^2}{2n}-\frac{\lambda^3}{3n^2}\Big)=\frac{e^{-\lambda} \lambda^4}{8 n^2}+o(n^{-2}), \ \ \mbox{ as } \ n\to\infty.
 \]
 }
 \end{REM}

 \begin{exam}
 \label{exam.Wasserstein}
 {\rm
  The distance $d_2$
  is helpful in proving stronger results, even for the
  ordinary
  Poisson approximation. For example, the
  Gini-Kantorovich, or Wasserstein,
  or transportation distance (see Novak (2019)) is defined by
  \begin{eqnarray*}
   d_W(S_n,Z_\lambda)
    &:=&
   \inf \E|S_n-Z_\lambda|
   \\
    &=&
   \mbox{$
   \sup_h \left|\E h(S_n)-\E h(Z_\lambda)\right|
   $}
   \\
    &=&
   \mbox{$
   \sum_{m=1}^\infty \Big|\Pr(S_n\geq m)-\Pr(Z_\lambda\geq m)\Big|,
   $}
  \end{eqnarray*}
  where the infimum is taken over all couplings of $S_n$ and $Z_{\lambda}$
  and the supremum over all functions $h:\{0,1,\ldots\}\to\R$ with $|h(m+1)-h(m)|\leq 1$  for all $m\geq 1$.
  For $g_1,g_2\in{\cal F}_2$ we consider a distance very
  similar to $d_2$, namely,
  \[
  \widetilde{d}_2(g_1,g_2):=\sum_{m=1}^{\infty}
  \frac{2^{m-1}}{(m-1)!} \Big|\mu_m(g_1)-\mu_m(g_2)\Big|.
  \]
  Then, $d_W(S_n,Z_\lambda)\leq \widetilde{d}_2(S_n,Z_\lambda)$.
  Indeed, a straightforward computation, using Theorem \ref{theo.d2}(a),
  yields
  \begin{eqnarray*}
  d_W(S_n,Z_\lambda)
   &\leq&
  \sum_{m=1}^\infty \sum_{k=m}^\infty \big|f_n(k)-g_{\lambda}(k)\big|
  =
  \sum_{k=1}^\infty k\big|f_n(k)-g_{\lambda}(k)\big|
  \\
   &=&
  \sum_{k=1}^\infty \frac{1}{(k-1)!}\left|
  \sum_{m=k}^\infty\frac{(-1)^{m-k}}{(m-k)!}
  (\mu_m-\lambda^m)\right|
  \leq
  \sum_{k=1}^\infty \frac{1}{(k-1)!}
  \sum_{m=k}^\infty\frac{\left|\mu_m-\lambda^m\right|}{(m-k)!}
  \\
   &=&
  \sum_{m=1}^\infty \frac{\left|\mu_m-\lambda^m\right|}{(m-1)!}
  \sum_{k=1}^m\frac{(m-1)!}{(k-1)! (m-k)!}=
  \widetilde{d}_2(S_n,Z_\lambda).
  \end{eqnarray*}
  On the other hand, using the lower bound from Lemma \ref{lem.comparison}
  and the simple fact that $\mu_m\leq \lambda^m$ for all $m$
  (the proof is omitted) we obtain
  $0\leq \lambda^m-\mu_m\leq \frac{(m)_2}{2}\lambda_2\lambda^{m-2}$,
  \[
  d_2(S_n,Z_\lambda)=\frac{1}{2}\sum_{m=1}^\infty \frac{2^m}{m!}
  \left\{ \lambda^m-\mu_m \right\}
  \leq
  \frac{1}{2}\sum_{m=1}^\infty \frac{2^m}{m!}
  \frac{(m)_2}{2}\lambda_2\lambda^{m-2}=e^{2\lambda} \lambda_2,
  \]
  and similarly, $\widetilde{d}_2(S_n,Z_\lambda)\leq 2(1+\lambda)e^{2\lambda}\lambda_2$. In this way we
  have produced
  a one-line
  proof of the inequality
  $d_W(S_n,Z_\lambda)\leq 2(1+\lambda)e^{2\lambda}\sum_{i=1}^n p_i^2$.
  }
  \end{exam}

 \begin{REM}
 \label{rem.Chen}
 {\rm
 Chen (1974) proved a stronger version
 of Poisson convergence
 for $S_n$. He showed that if
 $\sum_i p_i=\lambda$  and  $\max_i\{p_i\}\to 0$
 (equivalently, $\sum_i p_i^2\to 0$)
 then, as $n\to\infty$
 \[
  \sum_{k=0}^\infty h(k) \big|f_n(k)-g_\lambda(k)\big|
  \to 0
 \]
 for any function $h\geq 0$ for which $\E h(Z_\lambda)<\infty$.
 This mode of convergence is very strong; see Wang (1991).
 It is worth to point out that the moment inequality
 of Lemma \ref{lem.comparison}  offers a one-line proof
 of a slightly weaker result, namely, under the restriction
 that $h\geq 0$ and $\E Z_\lambda^2 h(Z_\lambda)<\infty$
 (this weaker mode is strong enough for all practical
 purposes). Indeed, we have
 \begin{eqnarray*}
 \sum_{k=0}^\infty h(k) \big|f_n(k)-g_\lambda(k)\big|
  &=&
 \sum_{k=0}^\infty \frac{h(k)}{k!}
 \left|\sum_{m=k}^\infty \frac{(-1)^{m-k}}{(m-k)!}
 \big(\mu_m-\lambda^m\big)\right|
 \\
  &\leq&
 \sum_{k=0}^\infty \frac{h(k)}{k!}
 \sum_{m=k}^\infty \frac{\lambda^m-\mu_m}{(m-k)!}
 \\
  &\leq&
 \frac{\lambda_2}{2}\sum_{k=0}^\infty \frac{h(k)}{k!}
 \sum_{m=k}^\infty \frac{m(m-1)\lambda^{m-2}}{(m-k)!}
 \\
  &=&
 \left(
 \frac{e^{2\lambda}}{2\lambda^2}\E\left\{
 h(Z_\lambda)\left(Z_{\lambda}^2+(2\lambda-1)Z_{\lambda}+\lambda^2
 \right)\right\}\right) \sum_{i=1}^n p_i^2\to 0.
 \end{eqnarray*}
  }
 \end{REM}

 \begin{REM}
 \label{rem.exact.d2}
 {\rm
 From the covariance identity we have
 \[
 \E \left\{ P_j(Z_\lambda)(Z_\lambda)_m\right\}
 =
 \lambda^j \E \left\{\Delta^j (Z_\lambda)_m\right\}
 =
 \lambda^j  (m)_j\E \left\{(Z_\lambda)_{m-j}\right\}
 =\lambda^m (m)_j.
 \]
 Hence, one advantage for considering $d_2$
 in connection with $\phi_{\nu}$
 when investigating the present problem is the fact that the
 corrected Poisson distribution admits an exact factorial
 moment sequence,
 \[
 \mbox{$
  \mu_m(\phi_\nu)=\lambda^m\left(1-\sum_{j=2}^{2\nu-2}
  \gamma_j \  (m)_j\right).
  $}
 \]
 Hence, whenever $\mu_m\leq \mu_m(\phi_{\nu})$ for all $m$, or
 $\mu_m\geq \mu_m(\phi_{\nu})$ for all $m$ (see Lemmas
 \ref{lem.comparison}, \ref{lem.comparison2}), we can derive
 a simple exact formula for $d_2$ as follows:
 \[
 \mbox{$
  d_2(f_n,\phi_\nu)=\frac{1}{2}\left|
  \sum_{m=0}^\infty \frac{2^m}{m!}\Big(\mu_m-\mu_m(\phi_\nu)\Big)\right|
  =
  \frac{1}{2}\left|
  \E\left\{ 3^{S_n}\right\} -
  e^{2\lambda}\left(1-\sum_{j=2}^{2\nu-2}\gamma_j
  \
  (2\lambda)^j\right)\right|.
  $}
 \]
  Therefore, from the independence of the indicators $I_i$ we get
  \[
  \mbox{$
   d_2(f_n,\phi_\nu)=\frac{1}{2}\left|
   \prod_{i=1}^n (1+2 p_i)-
  e^{2\lambda}\left(1-\sum_{j=2}^{2\nu-2}\gamma_j
  \
  (2\lambda)^j\right)\right|.
  $}
  \]
  }
 \end{REM}

 \begin{exam}
 \label{exam.binomial.exact}
 {\rm
 It is possible to obtain very accurate closed-form
 bounds of any order for corrected Poisson approximations
 in the simple case of the ordinary binomial distribution, i.e.,
 under the setup of Example \ref{exam.binomial}
 ($n>\lambda$, $p_i=\lambda/n$).
 We shall present some elementary results without proofs.

 We have $\mu_m=(n)_m\lambda^m/n^m$ and
 \[
  (n)_m=n^m-A_1n^{m-1}+A_2n^{m-2}+\cdots+(-1)^{m-1} A_{m-1} \ n
 \]
 where $A_k=A_k(m-1)=\sum_{1\leq i_1<\cdots<i_k\leq m-1}\prod_{j=1}^k i_j$
 are the well-known unsigned Stirling numbers of the first kind,
 see, e.g., Sibuya (1988);
 the usual notation is $S(m,k)=(-1)^{m-k} A_{m-k}$
 (for the signed numbers) and
 $\Big[\begin{array}{c}m \vspace{-.7ex}\\ k\end{array}\Big]
 =(-1)^{m-k}S(m,k)=|S(m,k)|$ for the unsigned ones.
 It can be verified that for any fixed $\nu\geq 1$ and $n\geq 1$,
 \begin{eqnarray*}
 (n)_m= n^m-A_1n^{m-1}+A_2n^{m-2}+\cdots+(-1)^{\nu-1} A_{\nu-1}
 n^{m-\nu+1}+(-1)^\nu R_\nu n^{m-\nu},
 \end{eqnarray*}
 where the remainder $0\leq R_\nu\leq A_{\nu}$.
 Also, the first few values of $A$'s are
 \begin{eqnarray*}
 &&
 \mbox{$
  A_1={m\choose 2}, \ A_2={m\choose 3}\frac{m-1}{4}, \
  A_3={m\choose 4}\frac{(m)_2}{2}, \
  A_4={m\choose 5}\frac{15 m^3-30 m^2+5 m+2}{48},
  $}
  \\
  &&
  \mbox{$
  A_5={m\choose 6}\frac{m(m-1)(3m^2-7m+2)}{16}, \
  A_6={m\choose 7}\frac{63m^5-315m^4+315m^3+91 m^2-42m-16}{576},
  $}
 \end{eqnarray*}
 and in general, $A_k$ is a polynomial of degree $2k$ in $m$ containing
 the binomial factor ${m\choose k+1}$ (so that $A_k=0$ for $k\geq m$).
 Our strategy in choosing a suitable $\phi_\nu$
 ($\nu\geq 2$, fixed) is quite simple: Equate the factorial moments
 $\mu_m(\phi_{\nu})$ of order
 up to $m=\nu$ with the truncated part of $\mu_m$, ignoring the
 remainder $R_\nu$. In this way, we obtain
 the equations
 \[
 \mbox{$
 \lambda^m\left(
 1-\frac{A_1}{n}+\frac{A_2}{n^{2}}+\cdots+(-1)^{\nu-1}
 \frac{A_{\nu-1}}{n^{\nu-1}}
 \right)
 =\mu_m(\phi_\nu)=\lambda^m
 \left( 1- \sum_{j=2}^{2\nu-2} \gamma_{\nu}(j) (m)_j\right).
 $}
 \]
 After the obvious cancelation of $\lambda^m$,
 both hands are polynomials in $m$ of degree $2\nu-2$, and, by
 solving a relative simple $(2\nu-3)\times(2\nu-3)$ linear system
 (with a triangular matrix of coefficients) we
 obtain the constants $\gamma_j(\nu)$,
 needed for the construction of the suitable
 corrected Poisson approximation
 of order $\nu$,
 \[
 \mbox{$
 \phi_{\nu}(k)=e^{-\lambda}\frac{\lambda^k}{k!}
 \left(1-\sum_{j=2}^{2\nu -2}\gamma_j(\nu)P_j(k)\right).
 $}
 \]
 Some values of $\gamma_j(\nu)$ are shown in the following table.
 \vspace{1em}

 \noindent
 {\small
 \begin{tabular}{r||llllll}
  {\small $j$} & {\small $\nu=2$} &
  {\small $\nu=3$} &
  {\small $\nu=4$} &
  {\small $\nu=5$} &
  {\small$\nu=6$} &
  {\small $\nu=7$}
  \\
  \hline
 $2$
   & $\frac{1}{2n}$
   & $\frac{1}{2n}$
   & $\frac{1}{2n}$
   & $\frac{1}{2n}$
   & $\frac{1}{2n}$
   & $\frac{1}{2n}$
   \vspace{.5ex}
   \\
 $3$
 &
 & $\frac{-1}{3n^2}$
 & $\frac{-1}{3n^2}$
 & $\frac{-1}{3n^2}$
 & $\frac{-1}{3n^2}$
 & $\frac{-1}{3n^2}$
  \vspace{.5ex}
  \\
 $4$
 &
 & $\frac{-1}{8n^2}$
 & $\frac{-1}{8n^2} + \frac{1}{4n^3}$
 & $\frac{-1}{8n^2} + \frac{1}{4n^3}$
 & $\frac{-1}{8n^2} + \frac{1}{4n^3}$
 & $\frac{-1}{8n^2} + \frac{1}{4n^3}$
  \vspace{.5ex}
  \\
 $5$ & & &
 $ \frac{1}{6n^3}$
 &
 $ \frac{1}{6n^3} - \frac{1}{5n^4}$
 &
 $ \frac{1}{6n^3} - \frac{1}{5n^4}$
 &
 $ \frac{1}{6n^3} - \frac{1}{5n^4}$
  \vspace{.5ex}
  \\
 $6$ & & &
 $\frac{1}{48n^3}$
 &
 $\frac{1}{48n^3} - \frac{13}{72n^4}$
 &
 $\frac{1}{48n^3} - \frac{13}{72n^4} + \frac{1}{6n^5}$
 &
 $\frac{1}{48n^3} - \frac{13}{72n^4} + \frac{1}{6n^5}$
  \vspace{.5ex}
  \\
 $7$ & & & &
 $\frac{-1}{24n^4}$
 &
 $\frac{-1}{24n^4} + \frac{11}{60 n^5}$
 &
 $\frac{-1}{24n^4} + \frac{11}{60 n^5} - \frac{1}{7n^6}$
  \vspace{.5ex}
  \\
 $8$ & & & &
 $\frac{-1}{384 n^4}$
 &
 $\frac{-1}{384 n^4} + \frac{17}{388n^5}$
 &
 $\frac{-1}{384 n^4} + \frac{17}{388n^5} - \frac{29}{160 n^6}$
  \vspace{.5ex}
  \\
 $9$ & & & & &
 $\frac{1}{144 n^5}$
 &
 $\frac{1}{144 n^5} - \frac{59}{810 n^6}$
  \vspace{.5ex}
  \\
 $10$ & & & & &
 $\frac{1}{3840n^5}$
 &
 $\frac{1}{3840n^5} - \frac{7}{576 n^6}$
  \vspace{.5ex}
  \\
 $11$ & & & & & &
 $\frac{-1}{1152 n^6}$
  \vspace{.5ex}
  \\
 $12$
 & & & & & &
 $\frac{-1}{46080 n^6}$
 \end{tabular}
 }
 {\small {\bf Table 1.} \rm Values $\gamma_j(\nu)$ needed in the correction
 $\phi_{\nu}$ as coefficients of orthogonal polynomials}.
 \vspace{1em}

 \noindent
 Since the remainder $R_{\nu}$ is nonnegative and bounded by
 $A_{\nu}$ (a polynomial in $m$ of degree $2\nu$),
 it follows that
 \[
 0\leq (-1)^{\nu}\Big(\mu_m-\mu_m(\phi_\nu)\Big)
 \leq A_\nu \frac{\lambda^m}{n^\nu}, \ \ m=1,2,\ldots \ .
 \]
 Thus,
 we obtain the exact formula (cf.\ Remark
 \ref{rem.exact.d2})
 \[
 d_2(\mbox{Bin}(n,\lambda/n),\phi_\nu)=\frac{1}{2}
 \Big|\Big(1+\frac{2\lambda}{n}\Big)^n-
 e^{2\lambda}\Big(1-\sum_{j=2}^{2\nu-2}\gamma_j(\nu)
 (2\lambda)^j\Big)\Big|
 \]
 and, also, the inequality
 \[
 d_2(\mbox{Bin}(n,\lambda/n),\phi_\nu)\leq C_{\nu}(\lambda) n^{-\nu},
 \ \ \ \
 C_\nu(\lambda)=\frac{1}{2}\sum_{m=\nu}^{\infty}
 \frac{(2\lambda)^m}{m!} A_{\nu}(m-1).
 \]
 Certainly, these results hold also for $\nu=1$; then
 $\phi_1=g_{\lambda}$
 is the usual Poisson$(\lambda)$ distribution
 and the constant is given by $C_1(\lambda)=\lambda^2 e^{2\lambda}$.
 By using the recurrence relation of Stirling numbers, see, e.g.,
 Sibuya (1988), it can be verified that
 $C_\nu(\lambda)=\lambda^{\nu+1}e^{2\lambda}Q_{\nu-1}(\lambda)$
 where $Q_{\nu}$ is a polynomial of degree $\nu$, $Q_0=1$, and
 $Q_{\nu}$ satisfies the recurrent differential equation
 \[
  \lambda Q_{\nu}'(\lambda)+(\nu+2)Q_{\nu}(\lambda)=
  2\lambda Q_{\nu-1}'(\lambda)+2(\nu+1+2\lambda)Q_{\nu-1}(\lambda), \ \
  \nu=1,2,\ldots;
 \]
 alternatively, one can derive the sequence $Q_\nu$ from the
 integral recurrence
 \[
  Q_{\nu+1}(\lambda)=
  2 Q_{\nu}(\lambda)+2\int_0^1 y^{\nu+2} (2\lambda y-1)
  Q_{\nu}(\lambda y) dy, \ \
  \nu=0,1,\ldots \ .
 \]
  The few first polynomials are
 $Q_1=\frac{4}{3}+\lambda$,
 $Q_2=2+\frac{8}{3}\lambda+\frac{2}{3}\lambda^2$,
 $Q_3=\frac{16}{5}+\frac{52}{9}\lambda+\frac{8}{3}\lambda^2
 +\frac{1}{3}\lambda^3$,
 $Q_4=\frac{16}{3} +\frac{176}{15}\lambda
 + \frac{68}{9}\lambda^2 +\frac{16}{9} \lambda^3 +
 \frac{2}{15} \lambda^4$. As a particular application, suppose we
 wish to approximate the binomial distribution $\mbox{Bin}(n,\lambda/n)$
 with a precision of order $n^{-4}$ for large $n$; that is, $\nu=4$.
 The proposed corrected Poisson approximation depends on the constants
 $\gamma_j(4)$, $j=2,\ldots,6$, appeared in column $\nu=4$
 of Table 1, i.e.,
 \[
 \mbox{$
 \phi_4(k)=e^{-\lambda}\frac{\lambda^k}{k!}
 \left(1
 -\frac{1}{2n} P_2(k)
 +\frac{1}{3n^2} P_3(k)
 +\Big(\frac{1}{8n^2}-\frac{1}{4n^3}\Big) P_4(k)
 -\frac{1}{6n^3} P_5(k)
 -\frac{1}{48n^3} P_6(k)
  \right).
  $}
 \]
 Then,
 \[
 d_2\Big(\mbox{Bin}(n,\lambda/n),\phi_4\Big)\leq C_{4}(\lambda) n^{-4}=
 \lambda^5 \Big(\frac{16}{5}+\frac{52}{9}\lambda+\frac{8}{3}\lambda^2
 +\frac{1}{3}\lambda^3\Big)e^{2\lambda} n^{-4}.
 \]
 Note that the bounds for $\nu=2$ and $\nu=3$
 follow from Theorems \ref{theo.main} and
 \ref{theo.main2}, respectively, when applied for equal $p_i$, and
 the constants are exactly the same.
 }
 \end{exam}

 It would be desirable to obtain general results for higher
 order (than three)
 corrected Poisson approximations for $S_n$,
 when the $p_i$'s are
 unequal. A natural conjecture, suggested by Example
 \ref{exam.binomial.exact} and Theorems
 \ref{theo.main}, \ref{theo.main2},
 is that for every $\nu\geq 2$,
 there exist (unique) constants $\gamma_2,\ldots,\gamma_{2\nu-2}$,
 depending only on $\lambda,\lambda_2,\ldots,\lambda_{\nu}$,
 such that
 \[
 d_2(f_n,\phi_\nu)\leq 
 Q_{\nu-1}(\lambda)e^{2\lambda}\sum_{i=1}^n p_i^{\nu+1},
 \]
 where $Q_{\nu}$ are the polynomials defined
 in Example \ref{exam.binomial.exact}.

\end{document}